\documentclass[12pt,a4paper]{article}
\usepackage{amsmath,amsfonts,amssymb,amscd}

\def\gr{\operatorname{gr}}
\def\Ker{\operatorname{Ker}}
\def\Ad{\operatorname{Ad}}
\def\ad{\operatorname{ad}}

\def\id{\operatorname{id}}

\def\d{\operatorname{d}}

\def\pr{\operatorname{pr}}
\def\Im{\operatorname{Im}}
\def\pt{\operatorname{pt}}
\def\Hom{\operatorname{Hom}}

\def\Aut{\operatorname{Aut}}
\def\Lie{\operatorname{Lie}}

\newcounter{th}
\def\t{\refstepcounter{th}{\bf \noindent{Theorem} \arabic{th}. }}

\newcounter{le}
\def\l{\refstepcounter{le}{\bf \noindent{Proposition} \arabic{le}. }}

\newcounter{lem}
\def\lem{\refstepcounter{lem}{\bf \noindent{Lemma} \arabic{lem}. }}

\newcounter{de}
\def\de{\refstepcounter{de}{\bf \noindent{Definition} \arabic{de}. }}

\newcounter{ex}
\def\ex{\refstepcounter{ex}{\bf \noindent{Example} \arabic{ex}. }}

\begin{document}

\begin{center}
\Large{\bf The splitting problem for complex homogeneous
supermanifolds \footnote{ Supported by Max Planck Institute for
Mathematics Bonn and AFR-grant, University of Luxembourg.
 } }
\end{center}

\begin{center}
     E.G. Vishnyakova
\end{center}

\begin{abstract}
It is a classical result that any complex analytic Lie supergroup
$\mathcal{G}$ is split \cite{kosz}, that is its structure sheaf is isomorphic to
the structure sheaf of a certain vector bundle. However, there do exist
non-split complex analytic homogeneous supermanifolds.

We study the question how to find out whether a complex analytic homogeneous supermanifold is split or non-split.
Our main result is a description of left invariant gradings on a complex analytic homogeneous
supermanifold $\mathcal{G}/\mathcal{H}$ in the terms of $\mathcal{H}$-invariants. As a corollary to our investigations we get some simple sufficient
conditions for a complex analytic homogeneous supermanifold to be split in terms of Lie algebras.

\end{abstract}

\bigskip

\section{Introduction}

A supermanifold is called split if its structure sheaf is isomorphic
to the exterior power of a certain vector bundle. By Batchelor's
Theorem any real supermanifold is non-canonically split. However, this
is false in the complex analytic case. The property of a
supermanifold to be split is very important for several reasons. For
instance, in \cite{DW} it was shown that the moduli space of super
Riemann surfaces is not projected (and in particular is not split)
for genus $g\geqslant 5$. The physical meaning of this result is
that \cite{DW}: "certain approaches to superstring perturbation
theory that are very powerful in low orders have no close analog in
higher orders".
 Another problem, when the property
of a supermanifold to be split is very important, is the calculation
of the cohomology group with values in a vector bundle over a
supermanifold. In the split case we may use the well understood
tools of complex analytic geometry. In the general case, several
methods were suggested by Onishchik's school: spectral sequences,
see e.g. \cite{OV}. All these methods connect the cohomology group
with values in a vector bundle with the cohomology group with values
in the corresponding split vector bundle.

How do we determine whether a complex analytic supermanifold is split or
non-split? Let me describe here some results in this direction that were obtained by Green,
Koszul, Onishchik and Serov. In \cite{Green} Green described a moduli space with a marked
point such that any non-marked point corresponds to a non-split
supermanifold while the marked point corresponds to a split
one. His idea was used for instance in \cite{DW}.
 The calculation of the Green moduli space is a
difficult problem itself, and in many cases the method is difficult to apply.
Furthermore,  Onishchik and Serov \cite{onigl, oniosp, onipisp}
consi\-dered  grading derivations, which correspond to
$\mathbb{Z}$-gradings of the structure sheaf of a supermanifold. For
example, it was shown that almost all super-grassmannians do not
possess such derivations, i.e. their structure sheaves do
not possess any $\mathbb{Z}$-gradings. Hence, in particular, they are
non-split. The idea of grading derivations was independently used by Koszul.  In \cite{Kosz} the following statement was proved: if the
tangent bundle of a supermanifold $\mathcal{M}$ possesses a
(holomorphic) connection then $\mathcal{M}$ is split. (Koszul's proof works in real and complex analytic cases.)  In fact, it
was shown that we can assign a grading derivation to any
supermanifold with a connection and that this grading derivation is
induced by a $\mathbb{Z}$-grading of a vector bundle.

Assume that a complex analytic supermanifold $\mathcal{M}=(\mathcal{M}_0,\mathcal O_{\mathcal{M}})$ is split. By definition this
means that its structure sheaf $\mathcal{O}_{\mathcal{M}}$ is isomorphic
to $\bigwedge \mathcal{E}$, where $\mathcal{E}$ is a locally free sheaf on the complex analytic manifold $\mathcal{M}_0$.
The sheaf $\bigwedge \mathcal{E}$ is naturally $\mathbb{Z}$-graded and
the isomorphism $\mathcal{O}_{\mathcal{M}} \simeq \bigwedge \mathcal{E}$
induces the $\mathbb{Z}$-grading in $\mathcal{O}_{\mathcal{M}}$. We call such gradings {\it split}.
The main result of our paper is a description of those left invariant
split gradings on a homogeneous superspace $\mathcal{G}/\mathcal{H}$
which are compatible with split gradings on $\mathcal{G}$.
  We also give sufficient conditions for pairs $(\mathfrak{g}, \mathfrak{h})$, where $\mathfrak{g} = \Lie \mathcal{G}$ and $\mathfrak{h} = \Lie \mathcal{H}$, such that $\mathcal{G}/\mathcal{H}$ is split.

\bigskip

\noindent  {\bf Acknowledgment.} The author is grateful to
A.~Onishchik, V.~Serganova, P.~Teichner and R.~Donagi for their attention
to this work and anonymous referee for useful comments.

\section{Complex analytic supermanifolds. Main definitions.} We will use the word "supermanifold" in the sense of Berezin and
Leites, see \cite{BL}, \cite{ley} and \cite{Man} for details.  Throughout, we will be interested in
the complex analytic version of the theory.
Recall that a {\it complex analytic superdomain of dimension $n|m$} is a $\mathbb{Z}_2$-graded
ringed space
$$
\mathcal{U} = \Big(U, \mathcal{F}_U \otimes \bigwedge(m) \Big),
$$
where $\mathcal{F}_U$ is the sheaf of holomorphic functions on an open set $U\subset \mathbb{C}^n$ and
$ \bigwedge(m)$ is the exterior (or Grassmann) algebra with $m$ generators.
 A {\it complex analytic supermanifold} of dimension $n|m$ is a $\mathbb{Z}_2$-graded ringed space that
is locally isomorphic to a complex superdomain of dimension $n|m$.

Let $\mathcal{M} = (\mathcal{M}_0,{\mathcal
O}_{\mathcal{M}})$ be a complex analytic supermanifold
 and
$$
\mathcal{J}_{\mathcal{M}} = ({\mathcal
O}_{\mathcal{M}})_{\bar 1} + ({\mathcal
O}_{\mathcal{M}})_{\bar 1}^2
$$
be the subsheaf of ideals generated
by odd
elements in ${\mathcal O}_{\mathcal{M}}$. We put $\mathcal{F}_{\mathcal{M}}:= {\mathcal
O}_{\mathcal{M}}/\mathcal{J}_{\mathcal{M}}$. Then $(\mathcal{M}_0,
\mathcal{F}_{\mathcal{M}})$ is a usual complex analytic manifold.  It
is called the {\it reduction} or {\it underlying space} of $\mathcal{M}$.
 We will
write $\mathcal{M}_0$ instead of $(\mathcal{M}_0,
\mathcal{F}_{\mathcal{M}})$ for simplicity of notation. Morphisms of
supermanifolds are just morphisms of the corresponding
$\mathbb{Z}_{2}$-graded ringed spaces. If $f: \mathcal{M} \to
\mathcal{N}$ is a morphism of supermanifolds, then we denote by
$f_0$ the
 morphism of the underlying spaces $\mathcal{M}_0 \to
\mathcal{N}_0$ and by $f^*$ the  morphism of the structure sheaves
$\mathcal{O}_{\mathcal{N}} \to (f_0)_*(\mathcal{O}_{\mathcal{M}})$.
 If $x\in \mathcal{M}_0$ and $\mathfrak m_x$ is
the maximal ideal of the local superalgebra $(\mathcal O_{\mathcal{M}})_x$, then
the vector superspace $T_x(\mathcal{M}):=(\mathfrak m_x/\mathfrak
m_x^2)^*$ is the tangent space of $\mathcal{M}$ at $x\in \mathcal{M}_0$.

Denote by ${\mathcal T}_{\mathcal{M}}$
the {\it tangent sheaf} or the {\it sheaf of vector fields} of
$\mathcal{M}$. In other words, ${\mathcal T}_{\mathcal{M}}$ is the
sheaf of derivations of the structure sheaf ${\mathcal
O}_{\mathcal{M}}$. Since the sheaf ${\mathcal O}_{\mathcal{M}}$ is
$\mathbb{Z}_2$-graded, the tangent sheaf ${\mathcal
T}_{\mathcal{M}}$ is also $\mathbb{Z}_2$-graded,
i.e. there is the natural decomposition ${\mathcal T}_{\mathcal{M}}
= ({\mathcal T}_{\mathcal{M}})_{\bar 0} \oplus ({\mathcal
T}_{\mathcal{M}})_{\bar 1}$, where
$$
({\mathcal T}_{\mathcal{M}})_{\bar i} :=  \Bigl\{ v\in {\mathcal T}_{\mathcal{M}}\,\, \big| \,
v\big(({\mathcal O}_{\mathcal{M}})_{\bar j}\big) \subset ({\mathcal O}_{\mathcal{M}})_{\bar j + \bar i}\,\Bigr\}.
$$

Let $\mathcal{M}_0$ be a complex analytic manifold and let $\mathcal{E}$ be the sheaf of holomorphic
sections of a vector bundle over $\mathcal{M}_0$. Then the ringed space
$(\mathcal{M}_0,\bigwedge\mathcal{E})$ is a supermanifold. In this case $\dim\,\mathcal{M} = n|m$, where $n =
\dim \mathcal{M}_0$ and $m$ is the rank of the locally free sheaf $\mathcal{E}$.

\medskip
\de\label{de split grading}
 A supermanifold $(\mathcal{M}_0,{\mathcal
O}_{\mathcal{M}})$ is called {\it split} if
${\mathcal
O}_{\mathcal{M}}\simeq \bigwedge\mathcal{E}$ for a
 locally free sheaf $\mathcal{E}$ on $\mathcal{M}_0$.
  The grading of $\mathcal{O}_{\mathcal{M}}$ induces by an isomorphism
  ${\mathcal O}_{\mathcal{M}}\simeq \bigwedge\mathcal{E}$ and  the natural $\mathbb{Z}$-grading
  of $\bigwedge \mathcal{E} = \bigoplus\limits_p \bigwedge\limits^p \mathcal{E}$ is called {\it split
   grading}.

\medskip

For example, all smooth supermanifolds are split by Batchelor's Theorem.
In \cite{Kosz} it was shown that all complex analytic Lie supergroups are split too.
 In this paper we study the splitting problem for complex analytic homogeneous supermanifolds.

\section{Lie supergroups and their homogeneous\\
 spaces}

\subsection{Lie supergroups and super Harish-Chandra pairs.} A
{\it Lie supergroup} is a group object in the category of supermanifolds,
i.e. it is a supermanifold $\mathcal{G}$ with three morphisms: the
multiplication morphism, the inversion morphism and the identity
morphism, which satisfy the usual conditions, modeling the group
axioms. In this case the underlying space $\mathcal{G}_0$ is a Lie group.  The structure sheaf of a (complex analytic) Lie
supergroup can be explicitly described in
 terms of the corresponding Lie superalgebra and underlying Lie group using
super Harish-Chandra pairs (see \cite{kosz} and \cite{V_TG} for more details). Let us describe this construction briefly.

\medskip

\de \label{de super Harish-Chandra pair} A {\it super Harish-Chandra pair} is a pair
$(\mathcal{G}_0,\mathfrak{g})$ that consists of a Lie group
$\mathcal{G}_0$ and a Lie superalgebra
$\mathfrak{g}=\mathfrak{g}_{\bar 0}\oplus\mathfrak{g}_{\bar 1}$ such
that $\mathfrak{g}_{\bar 0}= \Lie \mathcal{G}_0$ provided with a
representation $\Ad: \mathcal{G}_0\to \Aut \mathfrak{g}$ of
$\mathcal{G}_0$ in $\mathfrak{g}$ such that:
\begin{itemize}
\item $\Ad$ preserves
the parity and induces the adjoint representation of $\mathcal{G}_0$
on $\mathfrak{g}_{\bar 0}$;
\item the differential $(\d\Ad)_e$ at the identity $e\in \mathcal{G}_0$ coincides with
the adjoint representation $\ad$ of $\mathfrak g_{\bar 0}$ on $\mathfrak g$.
\end{itemize}

If a super Harish-Chandra pair $(\mathcal{G}_0,\mathfrak{g})$ is
given, it determines the Lie supergroup $\mathcal{G}$ in the
following way, see \cite{kosz}. Let $\mathfrak{U}(\mathfrak{g})$ be
the universal enveloping superalgebra of $\mathfrak{g}$. It is clear
that $\mathfrak{U}(\mathfrak{g})$ is a
$\mathfrak{U}(\mathfrak{g}_{\bar 0})$-module, where
$\mathfrak{U}(\mathfrak{g}_{\bar 0})$ is the universal enveloping
algebra of $\mathfrak{g}_{\bar 0}$. Recall that we denote by
$\mathcal{F}_{\mathcal{G}_0}$ the structure sheaf of the manifold
$\mathcal{G}_0$. The natural action of $\mathfrak{g}_{\bar 0}$ on
the sheaf $\mathcal{F}_{\mathcal{G}_0}$ gives rise to a structure of
$\mathfrak{U}(\mathfrak{g}_{\bar 0})$-module on
$\mathcal{F}_{\mathcal{G}_0}(U)$ for any open set $U\subset
\mathcal{G}_0$. Putting
$$
\mathcal{O}_{\mathcal{G}}(U) = \Hom_{\mathfrak{U}(\mathfrak{g}_{\bar
0})}\!\big(\mathfrak{U}(\mathfrak{g}), \mathcal{F}_{\mathcal{G}_0}(U)\big)
$$
for every open $U\subset \mathcal{G}_0$, we get a sheaf
$\mathcal{O}_{\mathcal{G}}$ of $\mathbb{Z}_2$-graded vector spaces. (Here we assume that the functions in
$\mathcal{F}_{\mathcal{G}_0}(U)$ are even.) The enveloping
superalgebra $\mathfrak{U}(\mathfrak{g})$ has a Hopf superalgebra
structure. Using this structure we can define the product of
elements from $\mathcal{O}_{\mathcal{G}}$ such that
$\mathcal{O}_{\mathcal{G}}$ becomes a sheaf of superalgebras, see
\cite{kosz} and \cite{V_TG} for details. A supermanifold structure on
$\mathcal{O}_{\mathcal{G}}$ is determined by the isomorphism
$\Phi_{\mathfrak{g}}:\mathcal{O}_{\mathcal{G}} \to \Hom \big(\bigwedge
(\mathfrak{g}_{\bar 1}), \mathcal{F}_{\mathcal{G}_0}\big)$, $f\mapsto
f\circ \gamma_{\mathfrak{g}}$, where
\begin{equation}
\label{isomorphism}
 \gamma_{\mathfrak{g}}: \bigwedge(\mathfrak{g}_{\bar 1})\to
\mathfrak{U}(\mathfrak{g}),\quad X_1\wedge \cdots \wedge X_r\mapsto
\frac{1}{r!}\sum_{\sigma\in S_r}(-1)^{|\sigma|} X_{\sigma(1)}\cdots
X_{\sigma(r)}.
\end{equation}
 The following formulas define the
multiplication morphism, the inversion morphism and the identity
morphism respectively:
\begin{equation}\label{umnozh}
\begin{split}
\mu^*(f)\left(X\otimes Y\right)(g,h)&=f\big(\Ad(h^{-1})(X)\cdot Y\big)\left(gh\right);\\
\iota^*(f)(X)(g)&=f\bigl(\Ad(g)(S(X))\bigr)\left(g^{-1}\right);\\
\varepsilon^*(f)&=f(1)(e).
\end{split}
\end{equation}
Here $X,Y\in\mathfrak{U}(\mathfrak{g}),\, f\in
\mathcal{O}_{\mathcal{G}},\, g,h\in \mathcal{G}_0$ and $S$ is the
antipode map of the Hopf superalgebra $\mathfrak{U}(\mathfrak{g})$.
Here we identify the enveloping superalgebra
$\mathfrak{U}(\mathfrak{g}\oplus\mathfrak{g})$ with the tensor
product $\mathfrak{U}(\mathfrak{g})\otimes
\mathfrak{U}(\mathfrak{g})$.

Sometimes we will identify the Lie superalgebra $\mathfrak{g}$ of a Lie
supergroup $\mathcal{G}$ with the tangent space $T_e(\mathcal{G})$
at $e\in \mathcal{G}_0$. The corresponding to $T\in T_e(\mathcal{G})$
left invariant vector field on $\mathcal{G}$ is given by
\begin{equation}\label{left invar vector field}
\left(\id \otimes T\right)\circ \mu^*,
\end{equation}
where $\mu$ is the multiplication morphism of $\mathcal{G}$.
(Recall that a vector field $Y$
on $\mathcal{G}$ is called {\it left invariant} if $(\id \otimes Y)\circ \mu^*= \mu^* \circ Y$.)
 Denote by $l_g$ and by $r_g$ the left and right
translations with respect to $g\in \mathcal{G}_0$, respectively.  The morphisms
$l_g$ and $r_g$ are given by the following formulas:
\begin{equation}\label{left right translations}
\begin{array}{l}
l_g^*(f)(X)(h) = f(X)(gh);\,\,\,\,\,r_g^*(f)(X)(h) =
f\big(\!\Ad(g^{-1})X\big)(hg),
\end{array}
\end{equation}
where $f\in \mathcal{O}_{\mathcal{G}}$, $X\in
\mathfrak{U}(\mathfrak{g})$ and $g,h\in \mathcal{G}_0$.

\subsection{Homogeneous supermanifolds.}

An {\it action of a Lie supergroup
$\mathcal{G}$ on a supermanifold} $\mathcal{M}$ is a
morphism $\nu:\mathcal{G}\times \mathcal{M}\to \mathcal{M}$ such that the usual conditions modeling group action axioms hold.
Any vector $X\in T_e(\mathcal{G})$ defines the vector field on $\mathcal{M}$ by the following formula:
\begin{equation}\label{fundamental vector field}
 X\mapsto
(X\otimes\, \id)\circ \nu^*.
\end{equation}

\smallskip

\de\label{def homogen supermnf}   An action $\nu$ is called {\it
transitive} if $\nu_{0}$ is a transitive action of the Lie group $\mathcal{G}_0$
 on $\mathcal{M}_0$ and the vector fields (\ref{fundamental vector field})
 generates the tangent space $T_x(\mathcal{M})$ at any point $x\in \mathcal{M}_0$.
 In this case the supermanifold
$\mathcal{M}$ is called $\mathcal{G}$-{\it homogeneous}.
 A supermanifold $\mathcal{M}$
is called {\it homogeneous}, if it possesses a transitive action of
a Lie supergroup.

\medskip

 If a supermanifold $\mathcal{M}$ is $\mathcal{G}$-homogeneous and
$\nu :\mathcal{G}\times \mathcal{M} \to \mathcal{M}$ is the
corresponding transitive action, then $\mathcal{M}$ is isomorphic to
the supermanifold $\mathcal{G}/\mathcal{H}$, where $\mathcal{H}$ is
the isotropy subsupergroup of a certain point (see \cite{V_Struc HSMF} for details). Recall that the
underlying space of $\mathcal{G}/\mathcal{H}$ is the complex analytic manifold
$\mathcal{G}_{0}/\mathcal{H}_0$ and the structure sheaf
$\mathcal{O}_{\mathcal{G}/\mathcal{H}}$ of $\mathcal{G}/\mathcal{H}$
is given by
\begin{equation}\label{homogeneous sheaf}
\mathcal{O}_{\mathcal{G}/\mathcal{H}} = \Bigl\{f\in
(\pi_0)_*(\mathcal{O}_{\mathcal{G}})\, \big|\, \mu_{\mathcal{G}\times
\mathcal{H}}^*(f) = \pr^*(f)\Bigr\},
\end{equation}
where $\pi_0:\mathcal{G}_0\to \mathcal{G}_0/\mathcal{H}_0$ is the
natural map, $\mu_{\mathcal{G}\times \mathcal{H}}$ is the
restriction of the multiplication map on $\mathcal{G}\times
\mathcal{H}$ and $\pr:\mathcal{G}\times \mathcal{H}\to \mathcal{G}$
is the natural projection. Using (\ref{umnozh}) we can rewrite the condition $\mu_{\mathcal{G}\times
\mathcal{H}}^*(f) = \pr^*(f)$  in the following way:
\begin{equation}\label{eq invariant condition}
 f\Big(\!\Ad\big(h^{-1}\big)(X)Y\Big)(gh) =
\left\{
                           \begin{array}{ll}
                             f(X)(g), & \hbox{$Y\in \mathbb{C}$;} \\
                             0, & \hbox{$Y \notin \mathbb{C}$;}
                           \end{array}
                         \right.
\end{equation}
where $X\in \mathfrak{U}(\mathfrak{g})$, $Y\in
\mathfrak{U}(\mathfrak{h})$, $\mathfrak{h} = \Lie \mathcal{H}$,
$g\in \mathcal{G}_0$ and $h\in \mathcal{H}_0$.

Let $Y\in \mathfrak{g}$ and $f\in \mathcal{O}_{\mathcal{G}}$.  Then
the operator defined by the formula
\begin{equation}\label{left inv vector field}
Y(f)(X) = (-1)^{p(Y)} f(XY),
\end{equation}
 where $p(Y)$ is the parity of $Y$, is a left invariant vector field on $\mathcal{G}$.
 From (\ref{left right translations}), (\ref{eq invariant condition}) and (\ref{left inv vector field}) it follows  that

 \smallskip

 {\it  $f\in
\mathcal{O}_{\mathcal{G}/\mathcal{H}}$ if and only if $f$ is
$\mathcal{H}_0$-right invariant, i.e. $r_h^*(f) = f$ for any $h\in
\mathcal{H}_0$, and $Y(f)=0$ for all $Y\in \mathfrak{h}_{\bar 1}$,
where $\mathfrak{h}= \mathfrak{h}_{\bar 0} \oplus \mathfrak{h}_{\bar
1}$.}

 \smallskip

Sometimes we will consider also the left action $\mathcal{H}\times
\mathcal{G}\to \mathcal{G}$ of a subsupergroup $\mathcal{H}$ on a
Lie supergroup $\mathcal{G}$. The corresponding quotient
supermanifold we will denote by $\mathcal{H}\backslash \mathcal{G}$.

\subsection{More about split supermanifolds.} Recall that a supermanifold
$\mathcal{M}$ is called  split if its structure sheaf
$\mathcal{O}_{\mathcal{M}}$ is isomorphic to $ \bigwedge \mathcal{E}$, where
$\mathcal{E}$ is a locally free sheaf on
$\mathcal{M}_0$. In this case, $\mathcal{O}_{\mathcal{M}}$ possesses
the $\mathbb{Z}$-grading induced by the natural $\mathbb{Z}$-grading
of $\bigwedge \mathcal{E} = \bigoplus\limits_p \bigwedge\limits^p \mathcal{E}$ and by isomorphism $\mathcal{O}_{\mathcal{M}} \simeq \bigwedge \mathcal{E}$.
Such gradings of $\mathcal{O}_{\mathcal{M}}$ we call split.

\medskip

\l\label{Ex_supergroups are nin-split} {\it Any Lie supergroup
$\mathcal{G}$ is split}.

\smallskip

This statement follows from the fact that any Lie supergroup is
determined by its super Harish-Chandra pair. A different proof of
this result (probably the first one) was given in \cite{Kosz}. For completeness we give here another proof.

\smallskip

\noindent {\it Proof.}  The underlying space
$\mathcal{G}_0$ is a closed Lie subsupergroup of $\mathcal{G}$. Hence,
there exists the homogeneous space $\mathcal{G}/\mathcal{G}_0$,
which is isomorphic to the supermanifold $\mathcal{N}$ such that
$\mathcal{N}_0$ is a point $\pt = \mathcal{G}_0/\mathcal{G}_0$ and
$\mathcal{O}_{\mathcal{N}}\simeq \bigwedge(m)$, where $m=\dim
\mathfrak{g}_{\bar 1}$. By definition, the structure sheaf
$\mathcal{O}_{\mathcal{N}}$ consists of all $r_g$-invariant
functions, $g\in \mathcal{G}_0$. We have the natural map
$\varphi: \mathcal{G}\to \mathcal{G}/\mathcal{G}_0$, where
$\varphi_{0}:\mathcal{G}_0 \to \pt$ and $\varphi^*: \mathcal{O}_{\mathcal{N}}
\to (\varphi_0)_*(\mathcal{O}_{\mathcal{G}})$ is the inclusion. It is
known that $\varphi: \mathcal{G}\to \mathcal{G}/\mathcal{G}_0$ is a
principal bundle (see \cite{V_Struc HSMF}). Using the fact that the underlying space of
$\mathcal{G}/\mathcal{G}_0$ is a point we get $\mathcal{G} \simeq
\mathcal{N}\times \mathcal{G}_0$. Note that this is an isomorphism of supermanifolds but not of Lie supergroups.$\Box$

\medskip

\ex\label{Ex_super-Grassmanians} As an example of a homogeneous
non-split supermanifold we can cite the super-grassmannian
$\mathbf{Gr}_{m|n,r|s}$ for $0<r<m$ and $0<s<n$. Super-grassmannians of other types are split (see Example \ref{ex split super-grassmanians}).

\medskip

Denote by $\verb"SSM"$ the category of split supermanifolds. Objects
$\operatorname{Ob}\, \verb"SSM"$ in this category are all split
supermanifolds $\mathcal{M}$ with fixed split gradings. Further if
$X,Y\in \operatorname{Ob}\, \verb"SSM"$, we put
$$
\Hom(X,Y)=  \Big\{\!
\begin{array}{l}
\text{morphisms from $X$ to $Y$}\\
\text{preserving the split
gradings}
\end{array}\!
\Big\}
$$

As in the category of supermanifolds, we can define in $\verb"SSM"$
a group object ({\it split Lie supergroup}), an action  of a split
Lie supergroup on a split supermanifold ({\it split action}) and a
 {\it split homogeneous supermanifold}.

There is a functor $\gr$ from the category of supermanifolds to the
category of split supermanifolds. Let us briefly describe this
construction. Let $\mathcal{M}$ be a supermanifold. Denote by
$\mathcal{J}_{\mathcal{M}}\subset \mathcal{O}_{\mathcal{M}}$ the
subsheaf of ideals generated by odd elements of
$\mathcal{O}_{\mathcal{M}}$. Then by definition $\gr
\mathcal{M}=(\mathcal{M}_0, \gr\mathcal{O}_{\mathcal{M}})$ is the
split supermanifold with the structure sheaf
$$
\gr\mathcal{O}_{\mathcal{M}}= \bigoplus_{p\geq 0}
(\gr\mathcal{O}_{\mathcal{M}})_p,\quad
\mathcal{J}_{\mathcal{M}}^0:=\mathcal{O}_{\mathcal{M}}, \quad
(\gr\mathcal{O}_{\mathcal{M}})_p:=
\mathcal{J}_{\mathcal{M}}^p/\mathcal{J}_{\mathcal{M}}^{p+1}.
$$
In this case $(\gr\mathcal{O}_{\mathcal{M}})_1$ is a locally free
sheaf and there is a natural isomorphism of
$\gr\mathcal{O}_{\mathcal{M}}$ onto $\bigwedge
(\gr\mathcal{O}_{\mathcal{M}})_1$. If
$\psi=(\psi_{0},\psi^*):\mathcal{M}\to \mathcal{N}$ is a morphism,
then $\gr(\psi)=(\psi_{0},\gr(\psi^*)) : \gr\mathcal{M}\to
\gr\mathcal{N}$ is defined by
$$
\gr(\psi^*)(f+\mathcal{J}_{\mathcal{N}}^p): =
\psi^*(f)+\mathcal{J}_{\mathcal{M}}^p \,\,\text{for}\,\, f\in
(\mathcal{J}_{\mathcal{N}})^{p-1}.
$$
Recall that by definition every morphism of supermanifolds is even
and as a consequence sends $\mathcal{J}_{\mathcal{N}}^p$ into
$\mathcal{J}_{\mathcal{M}}^p$.

\subsection{Split Lie supergroups.} Let $\mathcal{G}$ be a Lie
supergroup with the supergroup morphisms $\mu$, $\iota$ and
$\varepsilon$: the multiplication, the inversion and the identity morphism,
respectively. In this section we assign three split Lie supergroups $\mathcal{G}^1$,
$\mathcal{G}^2$ and $\mathcal{G}^3$ to $\mathcal{G}$ and we show that these split Lie supergroups are pairwise isomorphic.

{\bf (1)} The construction of $\mathcal{G}^1$ is very simple: we just apply functor $\gr$ to $\mathcal{G}$. Clearly,
 $\mathcal{G}^1: = \gr \mathcal{G}$ is a split Lie supergroup with
the supergroup morphisms $\gr(\mu)$, $\gr(\iota)$ and $\gr(\varepsilon)$.

\medskip

{\bf (2)} Consider the super
Harish-Chandra pair $(\mathcal{G}_0, \mathfrak{g}^2)$, where
$\mathfrak{g}^2$ is the following Lie superalgebra:
 $\mathfrak{g}^2$ and $\mathfrak{g}$ are isomorphic as vector superspaces and
the Lie bracket in $\mathfrak{g}^2$ is defined by the following formula:
\begin{equation}\label{bracket split}
[X,Y]_{\mathfrak{g}^2}=\left\{\!
        \begin{array}{ll}
          [X,Y]_{\mathfrak{g}}, & \hbox{if $X,Y\in \mathfrak{g}_{\bar 0}$ or $X\in \mathfrak{g}_{\bar 0}$ and $Y\in \mathfrak{g}_{\bar 1}$;} \\
          0, & \hbox{if $X,Y\in \mathfrak{g}_{\bar 1}$.}
        \end{array}
      \right.
\end{equation}
Denote by $\mathcal{G}^2$ the Lie supergroup corresponding to $(\mathcal{G}_0, \mathfrak{g}^2)$.

\medskip

{\bf (3)} Consider the sheaf $ \mathcal{O}_{\mathcal{G}^3} := \Hom_{\mathbb{C}}
(\bigwedge \mathfrak{g}_{\bar 1}, \mathcal{F}_{\mathcal{G}_0})$. For the ringed space $\mathcal{G}^3: = (\mathcal{G}_0, \mathcal{O}_{\mathcal{G}^3})$ we can repeat the construction from Section $3.1$. Indeed, this ringed space is clearly a supermanifold. Futhermore,  the exterior algebra $\bigwedge \mathfrak{g}_{\bar 1}$ is also a Hopf algebra. Therefore, we can define on $\mathcal{G}^3$  the multiplication, the inversion and the identity morphisms respectively by the following formulas:
\begin{equation}\label{umnozh_split}
\begin{split}
(\mu^3)^*(f)\big(X\wedge Y\big)(g,h)&=f\big(\Ad(h^{-1})(X)\wedge Y\big)(gh);\\
(\iota^3)^*\big(f\big)(X)(g)&=f\big(\Ad(g)(S'(X))\big)\big(g^{-1}\big);\\
(\varepsilon^3)^*\big(f\big)&=f\big(1\big)(e).
\end{split}
\end{equation}
Here $X,Y\in \bigwedge\mathfrak{g}_{\bar 1}$, $f\in \Hom_{\mathbb{C}}
(\bigwedge \mathfrak{g}_{\bar 1}, \mathcal{F}_{\mathcal{G}_0})$,
$g,\,h\in \mathcal{G}_0$ and $S'$ is the antipode map of the Hopf superalgebra
$\bigwedge\mathfrak{g}_{\bar 1}$. Hence, $\mathcal{G}^3 := (\mathcal{G}_0, \mathcal{O}_{\mathcal{G}^3})$ is a Lie supergroup.
Since
$$
\Hom_{\mathbb{C}}
(\bigwedge \mathfrak{g}_{\bar 1}, \mathcal{F}_{\mathcal{G}_0}) = \bigoplus\limits_{p\geq 0} \Hom_{\mathbb{C}}
(\bigwedge\limits^p \mathfrak{g}_{\bar 1}, \mathcal{F}_{\mathcal{G}_0})
$$
 is $\mathbb{Z}$-graded and the morphisms (\ref{umnozh_split})  preserve this $\mathbb{Z}$-grading, we see that $\mathcal{G}^3$ is a split Lie supergroup.

Later on we will need the explicit expression of left and right
translations $l'_g$ and $r'_g$ in $\mathcal{G}^3$:
\begin{equation}\label{left right translations split}
\begin{array}{l}
(l'_g)^*(f)(X)(h) = f(X)(gh);\,\,\,\,\,(r'_g)^*(f)(X)(h) =
f\big(\!\Ad(g^{-1})X\big)(hg),
\end{array}
\end{equation}
where $f\in \Hom_{\mathbb{C}}
(\bigwedge \mathfrak{g}_{\bar 1}, \mathcal{F}_{\mathcal{G}_0})$, $X\in
\bigwedge \mathfrak{g}_{\bar 1}$ and $g,h\in \mathcal{G}_0$.

\medskip

In fact, all these split Lie supergroups are isomorphic. To show this we need the following lemma:

\medskip

\lem\label{Lem_X,Y commute} {\it Let $\mathfrak{k}$ be a Lie
superalgebra and $X_i,Y_j\in \mathfrak{k}_{\bar 1}$, $i=1,\ldots,r$,
$j=1,\ldots,s$ be any elements. Assume that $[X_i,Y_j]=0$ for any $i,j$. Then we
have
$$
\gamma_{\mathfrak{k}}(X_1\wedge \cdots\wedge X_r\wedge Y_1\wedge
\cdots\wedge Y_s) = \gamma_{\mathfrak{k}}(X_1\wedge \cdots\wedge
X_r)\cdot \gamma_{\mathfrak{k}}(Y_1\wedge \cdots\wedge Y_s),
$$
where $\gamma_{\mathfrak{k}}$ is given by
$(\ref{isomorphism})$.

}

\medskip

\noindent {\it Proof.} A direct calculation.$\Box$

\medskip

\l\label{Prop mut of gr G} {\it  We have $\mathcal{G}^1 \simeq \mathcal{G}^2 \simeq \mathcal{G}^3$ in the category of Lie supergroups.}

\medskip

\noindent {\it Proof.} {\bf (a)} The statement $\mathcal{G}^1 \simeq \mathcal{G}^2$ was proven in \cite{V_Func}, Theorem $3$.

{\bf (b)} Let us show that $\mathcal{G}^2 \simeq \mathcal{G}^3$.
Applying Lemma \ref{Lem_X,Y commute} to $\mathfrak{g}^2$ and to any
elements $X_i,\,Y_j \in \mathfrak{g}^2_{\bar 1}$, we see that in this
case $\gamma_{ \mathfrak{g}^2}$ is not only isomorphism of super
coalgebras but of Hopf superalgebras. In other words, the isomorphism
$$
\Phi_{\mathfrak{g}^2} : \Hom_{\mathfrak{U}(\mathfrak{g}^2_{\bar 0})}
\big(\mathfrak{U}(\mathfrak{g}^2), \mathcal{F}_{\mathcal{G}_0}\big) \to \Hom_{\mathbb{C}}
\big(\bigwedge \mathfrak{g}_{\bar 1}, \mathcal{F}_{\mathcal{G}_0}\big)
$$
is an isomorphism of Lie supergroups.$\Box$

\section{Split grading operators}

Let again $\mathcal{M}$ be a supermanifold, $\gr\mathcal{M}$ be the
corresponding split supermanifold and $\mathcal{J}$ be the sheaf of
ideals generated by odd elements of $\mathcal{O}_{\mathcal{M}}$. We denote by
$\mathcal{T} = \mathcal{D}er \mathcal{O}_{\mathcal{M}}$ and by
$\gr{\mathcal{T}} = \mathcal{D}er (\mathcal{O}_{\gr\mathcal{M}})$
the tangent sheaf  of $\mathcal{M}$ and of $\gr\mathcal{M}$, respectively. The sheaf $\mathcal{T}$ is
naturally $\mathbb{Z}_2$-graded and the sheaf $\gr\mathcal{T}$ is
naturally $\mathbb{Z}$-graded: the gradings are induced by the
$\mathbb{Z}_2$ and $\mathbb{Z}$-grading of
$\mathcal{O}_{\mathcal{M}}$ and $\gr\mathcal{O}_{\mathcal{M}}$,
respectively. In other words, we have the decomposition:
$$
\mathcal{T} = \mathcal{T}_{\bar 0} \oplus \mathcal{T}_{\bar 1}, \quad \gr\mathcal{T} = \bigoplus\limits_{p\geq -1} (\gr\mathcal{T})_p.
$$

The sheaves $\mathcal{T}$ and $ \gr\mathcal{T}$ are related: this relation can be expressed by the following exact sequence:
\begin{equation}\label{form_tangent sheaf_gen}
0\longrightarrow \mathcal{T}_{(2)\bar 0} \longrightarrow \mathcal{T}_{\bar 0}
\stackrel{\alpha}{\longrightarrow} (\gr{\mathcal{T}})_{0} \longrightarrow 0,
\end{equation}
where
$$
\mathcal{T}_{(2)\bar 0} = \{ v\in \mathcal{T}_{\bar 0} \mid
v(\mathcal{O}_{\mathcal{M}}) \subset \mathcal{J}^2\}.
$$
The morphism $\alpha$ in
(\ref{form_tangent sheaf_gen}) is the composition of the natural morphism
$\mathcal{T}_{\bar 0} \to \mathcal{T}_{\bar 0}/\mathcal{T}_{(2)\bar
0}$ and the isomorphism $\mathcal{T}_{\bar 0}/\mathcal{T}_{(2)\bar 0} \to
(\gr{\mathcal{T}})_{0}$ that is given by
$$
[w] \longmapsto  \tilde{w},\quad
\tilde{w}\big(f+\mathcal{J}^{p+1}\big):= w(f) + \mathcal{J}^{p+1},
$$
where $w\in \mathcal{T}_{\bar 0}$, $[w]$ is the image of $w$ in
$\mathcal{T}_{\bar 0}/\mathcal{T}_{(2)\bar 0}$ and $f\in
\mathcal{J}^{p}$.

Assume that the sheaf $\mathcal{O}_{\mathcal{M}}$ is
$\mathbb{Z}$-graded, i.e. $\mathcal{O}_{\mathcal{M}} = \bigoplus\limits_p
(\mathcal{O}_{\mathcal{M}})_p$. Then we have the map $w:\mathcal{O}_{\mathcal{M}} \to \mathcal{O}_{\mathcal{M}}$
defined by $w(f) = pf$, where $f\in (\mathcal{O}_{\mathcal{M}})_p$. Such maps
are called {\it grading operators } on $\mathcal{M}$.

\medskip

\de\label{de split grading operator}We call a grading operator
$w$ on $\mathcal{M}$ a {\it split grading operator} if it
corresponds to a split grading of $\mathcal{O}_{\mathcal{M}}$, see Definition \ref{de split grading}.

\medskip

 In fact any split grading operator $w$ on $\mathcal{M}$ is
an even vector field on $\mathcal{M}$. Indeed, $w$ is linear, it preserves
the parity in $\mathcal{O}_{\mathcal{M}}$ and for $f\in (\mathcal{O}_{\mathcal{M}})_p$ and $g\in (\mathcal{O}_{\mathcal{M}})_q$ we have:
$$
w(fg) = (p+q) fg = (pf) g + f(qg) = w(f)g + fw(g).
$$
Note that $fg \in (\mathcal{O}_{\mathcal{M}})_{p+q}$.

By definition the sheaf $\gr\mathcal{O}_{\mathcal{M}}$ is $\mathbb{Z}$-graded. Denote by $a$ the corresponding split grading operator.

\medskip

\lem\label{Prop_krit_v lifts} {\it 1. A supermanifold $\mathcal{M}$
is split if and only if the vector field $a$ is contained in $\Im H^0(\alpha)$, where
$$
H^0(\alpha): H^0(\mathcal{M}_0,\mathcal{T}_{\bar 0}) \to
H^0(\mathcal{M}_0,(\gr{\mathcal{T}})_{0}).
$$
(We applied the functor $H^0(\mathcal{M}_0,-)$ to the sequence (\ref{form_tangent sheaf_gen}).
We write $H^0(\alpha)$ instead of $H^0(\mathcal{M}_0,\alpha)$ for notational simplicity.)

2. If $w$ is a split grading operator on $\mathcal{M}$, then any
other split grading operator on $\mathcal{M}$ has the form $w+
\chi$, where  $\chi \in H^0(\mathcal{M}_0,\mathcal{T}_{(2)\bar 0})$.

}
\medskip

 \noindent{\it Proof.} 1. The statement of the lemma can be
deduced from the following observation made by Koszul in \cite[Lemma
$1.1$ and Section $3$]{Kosz}. Let $A$ be a commutative superalgebra
over $\mathbb{C}$ and $\mathfrak{m}$ be a nilpotent ideal in $A$. An
even derivation $w$ of $A$ is called {\it adapted to the filtration}
$$ A\supset \mathfrak{m} \supset \mathfrak{m}^2 \ldots
$$
 if
$$
(w-r\id)(\mathfrak{m}^r) \subset \mathfrak{m}^{r+1} \,\,\,\text{for
any $r\geq 0$.}
$$
Denote by $D_{\mathfrak{m}}^{ad}$ the set of all derivations adapted
to $\mathfrak{m}$. In \cite[Lemma $1.1$]{Kosz} it was shown that
$D_{\mathfrak{m}}^{ad}$ is not empty if and only if the filtration of $A$ is
splittable. Moreover, if $w\in D_{\mathfrak{m}}^{ad}$, then the corresponding splitting of $A$ is given
by eigenspaces of the derivation $w$: $A= \bigoplus_i A_i$, where $A_i$ is
the eigenspace of $w$ with the eigenvalue $i$, and $\mathfrak{m}^r= A_r\oplus \mathfrak{m}^{r+1}$ for all $r\geq 0$.

We apply Koszul's observation to the sheaf of superalgebras $\mathcal{O}_{\mathcal{M}}$ and
its subsheaf of ideals $\mathcal{J}$. The set
$D_{\mathcal{J}}^{ad}$ is in this case the set of global derivations
of
$\mathcal{O}_{\mathcal{M}}$ adapted to the filtration
\begin{equation}\label{eq filtration in structure sheaf}
\mathcal{O}_{\mathcal{M}}\supset \mathcal{J} \supset \mathcal{J}^2 \supset \ldots.
\end{equation}

Clearly, $D_{\mathcal{J}}^{ad}$ is not empty if and only if $a$ is
contained in $\Im H^0(\alpha)$. (Actually,
$H^0(\alpha)(D_{\mathcal{J}}^{ad})=a$.) Furthermore, if the
supermanifold $\mathcal{M}$ is split, i.e. we have a split grading
$\mathcal{O}_{\mathcal{M}} = \bigoplus_{p\geq 0}
(\mathcal{O}_{\mathcal{M}})_p$, then $\mathcal{J}^q=
\bigoplus_{p\geq q} (\mathcal{O}_{\mathcal{M}})_p$ and
$\mathcal{J}^q = (\mathcal{O}_{\mathcal{M}})_q \oplus
\mathcal{J}^{q+1}$. Hence, the split grading determine the splitting
of the filtration (\ref{eq filtration in structure sheaf}) and the
corresponding split grading operator belongs to
$D_{\mathcal{J}}^{ad}$.

Conversely, if there exists $w\in D_{\mathcal{J}}^{ad}$, then we can
decompose the sheaf $\mathcal{O}_{\mathcal{M}}$ into eigenspaces
$$
(\mathcal{O}_{\mathcal{M}})_q:= \{f\in \mathcal{O}_{\mathcal{M}} | w(f) = qf \}.
$$
In this case the sheaves $\bigoplus\limits_p (\mathcal{O}_{\mathcal{M}})_p$ and $\gr \mathcal{O}_{\mathcal{M}}$ are isomorphic as
 $\mathbb{Z}$-graded sheaves of superalgebras since $\mathcal{J}^q = (\mathcal{O}_{\mathcal{M}})_q \oplus
 \mathcal{J}^{q+1}$. Hence, the supermanifold is split.

\smallskip

2.  Applying  the left-exact functor $H^0(\mathcal{M}_0,-)$ to (\ref{form_tangent sheaf_gen}), we get the following exact sequence:
$$
0\longrightarrow H^0(\mathcal{M}_0,\mathcal{T}_{(2)\bar 0}) \longrightarrow H^0(\mathcal{M}_0,\mathcal{T}_{\bar 0})
\stackrel{H^0(\alpha)}{\longrightarrow} H^0(\mathcal{M}_0,(\gr{\mathcal{T}})_{0}).
$$
If $w_1,\,w_2$ are two split grading operators on $\mathcal{M}$, then
$$
H^0(\alpha)(w_1) = H^0(\alpha)(w_2) = a,
$$
according to the part 1.
Therefore, $w_1-w_2 \in H^0(\mathcal{M}_0,\mathcal{T}_{(2)\bar 0})$.
 The
result follows.$\Box$

\medskip

\ex\label{Ex_split grading operator for split Lie sgroup} Consider
the supermanifold $\mathcal{G}_0\backslash \mathcal{G}$. Its
structure sheaf is isomorphic to $\bigwedge(\mathfrak{g}_{\bar 1})$
(compare with Example \ref{Ex_supergroups are nin-split}). Denote by
$(\varepsilon_i)$ the system of odd (global) coordinates on
$\mathcal{G}_0\backslash \mathcal{G}$. An example of a split grading
operator on the Lie supergroup $\mathcal{G}$ is $\sum \varepsilon^i
X_i$. Here $(X_i)$ is a basis of odd left invariant vector fields on
$\mathcal{G}$ such that $X_i(\varepsilon^j)(e) = \delta_i^j$. We may
produce other examples if we use right invariant vector fields or
odd (global) coordinates on $\mathcal{G}/ \mathcal{G}_0$.

By Lemma \ref{Prop_krit_v lifts}, any split grading operator on a
Lie supergroup $\mathcal{G}$ is given by $\sum \varepsilon^i X_i +
\chi$, where $\chi \in H^0(\mathcal{G}_0,\mathcal{T}_{(2)\bar 0})$
is any vector field on $\mathcal{G}$.

\section{Compatible split gradings on $\mathcal{G}/\mathcal{H}$}

\subsection{Compatible gradings on $\mathcal{G}/\mathcal{H}$.}

Let $\mathcal{G}$ be a Lie supergroup and $\mathcal{M} = \mathcal{G}/\mathcal{H}$ be a homogeneous supermanifold.
As above we denote by $\pi: \mathcal{G}\to  \mathcal{G}/\mathcal{H}$ the natural projection.

\medskip

\de\label{de compatible grading} A split grading of the sheaf
$\mathcal{O}_{\mathcal{G}} = \bigoplus\limits_p (\mathcal{O}_{\mathcal{G}})_p$ is
called {\it compatible} with the inclusion $\mathcal{O}_{\mathcal{M}}\subset
(\pi_0)_* (\mathcal{O}_{\mathcal{G}})$ if the following holds:
$$
f\in \mathcal{O}_{\mathcal{M}} \,\,\, \Rightarrow f_p\in
\mathcal{O}_{\mathcal{M}}\,\, \text{for all} \,\,\, p,
$$
where $f= \sum f_p$ and $f_p\in (\pi_0)_*((\mathcal{O}_{\mathcal{G}})_p)$.

\medskip

Let us take any split grading operator $w$ on $\mathcal{G}$. Clearly,
the corresponding split grading of $\mathcal{O}_{\mathcal{G}}$ is compatible with $\mathcal{O}_{\mathcal{M}}$ if and only if
$w(\mathcal{O}_{\mathcal{M}}) \subset \mathcal{O}_{\mathcal{M}}$.
It is not clear from Definition \ref{de compatible grading} that the compatible grading
\begin{equation}\label{eq compatible grading}
(\mathcal{O}_{\mathcal{M}})_p= \mathcal{O}_{\mathcal{M}}  \cap
(\pi_0)_*((\mathcal{O}_{\mathcal{G}})_p)
\end{equation}
of $\mathcal{O}_{\mathcal{M}}$, if it exists, is a
split grading of $\mathcal{O}_{\mathcal{M}}$. However, the following proposition holds:

\medskip

\l\label{if compatible then good} {\it Assume that we have the $\mathbb{Z}$-grading:
$$
\mathcal{O}_{\mathcal{M}} = \bigoplus\limits_{p\geq 0} (\mathcal{O}_{\mathcal{M}})_p,
$$
 where
$(\mathcal{O}_{\mathcal{M}})_p$ are as in (\ref{eq compatible grading}). Then this grading is a split grading.}

\medskip

\noindent{\it Proof.} The idea of the proof is to apply Lemma \ref{Prop_krit_v lifts} to
the grading operator $w':= w|_{\mathcal{O}_{\mathcal{M}}}$ on $\mathcal{M}$.
Denote by $\mathcal{J}_{\mathcal{M}}$ and by $\mathcal{J}_{\mathcal{G}}$ the sheaves of
ideals generated by odd elements of $\mathcal{O}_{\mathcal{M}}$ and $\mathcal{O}_{\mathcal{G}}$, respectively. Our aim is to show that
$$
w'(f)+ \mathcal{J}_{\mathcal{M}}^{p+1} = pf + \mathcal{J}_{\mathcal{M}}^{p+1},
$$
where $f\in\mathcal{J}_{\mathcal{M}}^{p}$. In other words, we want to show that $H^0(\alpha) (w')$
is a split grading operator for the grading of $\gr\mathcal{O}_{\mathcal{M}}$. (We use notations of Lemma \ref{Prop_krit_v lifts}.) We have:
$$
\begin{array}{c}
(\gr \pi)^*(w'(f) + \mathcal{J}_{\mathcal{M}}^{p+1}) = w(f) + \mathcal{J}_{\mathcal{G}}^{p+1} = pf + \mathcal{J}_{\mathcal{G}}^{p+1};\\
 \rule{0pt}{5mm}(\gr \pi)^*( pf + \mathcal{J}_{\mathcal{M}}^{p+1}) = pf + \mathcal{J}_{\mathcal{G}}^{p+1}.
\end{array}
$$
Since the map $(\gr \pi)^*$ is injective, we get, $w'(f)+
\mathcal{J}_{\mathcal{M}}^{p+1} = pf +
\mathcal{J}_{\mathcal{M}}^{p+1}$. $\Box$

\subsection{$\mathcal{H}$-invariant split grading operators.}

 First of all let us consider the situation when a split
grading operator $w$ on $\mathcal{G}$ is invariant with respect to a
Lie subsupergroup $\mathcal{H}$. In terms of super Harish-Chandra pairs this means:
\begin{equation}\label{v is H invariant}
  \begin{array}{ll}
    r_h^* \circ w = w \circ r_h^*, & \hbox{for all\,\,\, $h\in \mathcal{H}_0$;} \\

 [Y,w]=0, & \hbox{for all $Y\in \mathfrak{h}_{\bar 1}$.}
  \end{array}
\end{equation}
Here $(\mathcal{H}_0,\mathfrak{h})$ is the super Harish-Chandra pair
of $\mathcal{H}$, $r_h$ is the right translation and $Y$ is an odd left invariant vector field.

\medskip

\l\label{H inveriant grading operators} {\it Assume that  $w$
is an $\mathcal{H}$-invariant split grading
operator on $\mathcal{G}$, i.e. equations
(\ref{v is H invariant}) hold. Then
 $\mathcal{H}$ is an ordinary Lie group.}

\medskip

\noindent{\it Proof.} The idea of the proof is to show that the Lie
superalgebra $\mathfrak{h}$ of $\mathcal{H}$ has the trivial odd part: $\mathfrak{h}_{\bar 1}=\{0\}$.

In Example \ref{Ex_split grading operator for
split Lie sgroup} we saw that any split grading operator on $\mathcal{G}$ is
given by $w= \sum \varepsilon^i X_i + \chi$. If $Z$ is a vector
field on $\mathcal{G}$, denote by $Z_e\in T_e(\mathcal{G})$ the
corresponding tangent vector at the identity $e\in \mathcal{G}_0$. Consider the
second equation in (\ref{v is H invariant}). At the point $e$, we
have
$$
\begin{array}{c}
[Y,w]_e = \big(\!\sum\limits_i Y(\varepsilon^i) X_i - \sum\limits_i \varepsilon^i Y\circ X_i
- \sum\limits_i \varepsilon^i X_i\circ Y + [Y,\chi]\big)_e =0
\end{array}
$$
for any $Y\in \mathfrak{h}_{\bar 1}.$ Furthermore,
$$
(\sum\limits_i \varepsilon^i Y\circ X_i
- \sum\limits_i \varepsilon^i X_i\circ Y)_e =0 \quad\text{and} \quad [Y,\chi]_e =0,
$$
because $\varepsilon^i(e) =0$ and because $\chi\in  H^0(\mathcal{M}_0,\mathcal{T}_{(2)\bar 0})$. Therefore,
$$
[Y,w]_e =  \sum\limits_i
Y(\varepsilon^i)(e) (X_i)_e =0
$$
The tangent vectors $(X_i)_e$ form a basis in $T_e(\mathcal{G})_{\bar 1}$,
hence $Y(\varepsilon^i)(e)=0$ for all $i$. The last statement is equivalent
to $Y_e=0$. Since $Y$ is a left invariant vector field, we get $Y=0$. The proof is complete.$\Box$

\medskip

\noindent {\bf Remark.}  It is well known that the supermanifold
$\mathcal{G}/\mathcal{H}$, where $\mathcal{H}$ is an ordinary Lie group,
is split (see \cite{kosz} or \cite{V_TG}). Therefore, the case of
$\mathcal{H}$-invariant split grading operators does not lead to new examples of homogeneous split supermanifolds.

\subsection{$\mathcal{G}_{0}$-left invariant split grading operators.}

Consider now a more general situation, when a split grading operator
$w$ leaves $\mathcal{O}_{\mathcal{M}}$ invariant. Let $f\in
\mathcal{O}_{\mathcal{M}}$. Then $w(f)\in \mathcal{O}_{\mathcal{M}}$
if and only if
$$
r^*_h(w(f))= w(f)\quad \text{and} \quad Y(w(f))=0
$$
for $h\in \mathcal{H}_0$
and $Y\in \mathfrak{h}_{\bar 1}$. These conditions are equivalent to
the following ones:
\begin{equation}\label{conditions v invariant}
\begin{array}{c}
(r^*_h \circ w \circ (r^{-1}_h)^* - w )|_{\mathcal{O}_{\mathcal{M}}}
=0;\quad [Y,w]|_{\mathcal{O}_{\mathcal{M}}} = 0.
\end{array}
\end{equation}
Recall that $r^{-1}_h = r_{h^{-1}}$.

It seems to us that the system (\ref{conditions v invariant}) is hard
to solve in general. Consider now a special type of split grading operators, called $\mathcal{G}_0$-left invariant grading operators.

\medskip

\de A split grading of $\mathcal{O}_{\mathcal{G}}$ is
called {\it $\mathcal{G}_0$-left invariant} if it is invariant with
respect to left translations. In other words, from $f\in
(\mathcal{O}_{\mathcal{G}})_p$ it follows that $l^*_g(f)\in
(\mathcal{O}_{\mathcal{G}})_p$ for all $g\in \mathcal{G}_0$.

\medskip

It is easy to see that a split grading of
$\mathcal{O}_{\mathcal{G}}$ is $\mathcal{G}_0$-left invariant if and
only if the corresponding split grading operator $w$ is invariant
with respect to left translations: $l^*_g\circ w = w \circ l^*_g$,
$g\in \mathcal{G}_0$. For example, the split grading operator $\sum
\varepsilon^i X_i$ constructed in Example \ref{Ex_split grading
operator for split Lie sgroup} is a $\mathcal{G}_0$-left invariant
split grading operator, because $\varepsilon^i$ are
$\mathcal{G}_0$-left invariant functions and $X_i$ are left
invariant vector fields. In this section we will describe all such
operators.

In Section $3.4$ we have seen that the supermanifold
$(\mathcal{G}_0,\Hom_{\mathbb{C}}(\bigwedge \mathfrak{g}_{\bar 1}, \mathcal{F}_{\mathcal{G}_0}))$ is a Lie
supergroup isomorphic to $\gr \mathcal{G}$. We need the following lemma:

\medskip

\lem\label{lem Phi is left right invariant} {\it The map
$$
\begin{array}{rl}
\Phi_{\mathfrak{g}}: \mathcal{O}_{\mathcal{G}}& \to \Hom_{\mathbb{C}}(\bigwedge \mathfrak{g}_{\bar 1}, \mathcal{F}_{\mathcal{G}_0}),\\
f&\mapsto f\circ \gamma_{\mathfrak{g}}
\end{array}
$$
from Section $3.1$ is invariant with
respect to left and right translations.

 }

\medskip

\noindent{\it Proof.} For any $h\in
\mathcal{G}_0$, denote by $r'_h$ and $l'_h$ the right and the left translation in the Lie
supergroup $\mathcal{G}^3 =(\mathcal{G}_0,\Hom_{\mathbb{C}}(\bigwedge \mathfrak{g}_{\bar 1}, \mathcal{F}_{\mathcal{G}_0}))$,
respectively. (See, (\ref{left right translations split})) Let us show that
\begin{equation}\label{r_g Phi = Phi R_g}
(r'_h)^* \circ \Phi_{\mathfrak{g}} = \Phi_{\mathfrak{g}} \circ
r_h^*.
\end{equation}
Let us take $Z\in \bigwedge \mathfrak{g}_{\bar 1}$ and $g,h\in
\mathcal{G}_0$. Using (\ref{left right translations}) we have
$$
\begin{array}{c}
\big[(r'_h)^* \circ \Phi_{\mathfrak{g}}\big](f) (Z)(g) =
 \Phi_{\mathfrak{g}}(f)\big(\!\Ad(h^{-1})(Z)\big)(gh) =\\
 \rule{0pt}{5mm}f\big(\!\gamma_{\mathfrak{g}}(\Ad(h^{-1})(Z))\big)(gh) =
 f\big(\!\Ad(h^{-1})(\gamma_{\mathfrak{g}}(Z))\big)(gh) =\\
 \rule{0pt}{5mm}r_h^*(f)\big(\!\gamma_{\mathfrak{g}}(Z)\big)(g) = [\Phi_{\mathfrak{g}} \circ
r_h^*](f)(Z)(g).
\end{array}
$$
Similarly, we get
$$
(l'_h)^* \circ \Phi_{\mathfrak{g}} =
\Phi_{\mathfrak{g}} \circ l_h^*.
$$
\begin{flushright}
$\Box$
\end{flushright}

\medskip

The following observation is known to experts, but we cannot find it in the literature:

\medskip

\lem\label{left invariant grading operators} {\it The space of $\mathcal{G}_0$-left invariant
 vector fields $H^0(\mathcal{G}_0, \mathcal{T})^{\mathcal{G}_0}$ on a Lie supergroup $\mathcal{G}$
 is isomorphic to  $H^0(\pt, \mathcal{O}_{\mathcal{G}_0\backslash
\mathcal{G}}) \otimes \mathfrak{g}$. The isomorphism is given by:
$$
f\otimes Z \stackrel{F}{\longmapsto}  fZ,
$$
where $f\in H^0(\pt, \mathcal{O}_{\mathcal{G}_0\backslash
\mathcal{G}})$ and $Z \in \mathfrak{g}$.
 }

\medskip

\noindent{\it Proof.} Clearly, the map $F$ is injective and its image is contained in the
vector space $H^0(\mathcal{G}_0, \mathcal{T})^{\mathcal{G}_0}$.  Let us show that any vector
field $v$ in $H^0(\mathcal{G}_0, \mathcal{T})^{\mathcal{G}_0}$ is contained in $\Im(F)$.

Let $(X_i)$ and $(Z_j)$ be a basis of odd and even left invariant (with respect to the
supergroup $\mathcal{G}$) vector fields on $\mathcal{G}$, respectively. Assume that
$$
v= \sum f^iX_i + \sum g^j Z_j,
$$
where $f^i, g^j \in H^0(\mathcal{G}_0,\mathcal{O}_{\mathcal{G}})$,
be the decomposition of $v$ with respect to this basis.
 We
have:
$$
\begin{array}{rl}
l_g^*\circ v = &\sum
l_g^*(f^i) l_g^* \circ X_i + \sum l_g^*(g^j)  l_g^*\circ Z_j=
\\
&\rule{0pt}{5mm}\sum
l_g^*(f^i)  X_i\circ l_g^* + \sum l_g^*(g^j) Z_j \circ l_g^* =
v\circ l_g^*.
\end{array}
$$
Therefore,
$l_g^*(f^i) = f^i$ and $l_g^*(g^j) = g^j$ for all $g\in \mathcal{G}_0$. In other words, $f^i,g^j\in H^0(\pt, \mathcal{O}_{\mathcal{G}_0\backslash
\mathcal{G}})$.  The proof
is complete.$\Box$

\medskip

The Lie supergroup $\mathcal{G}$ acts on the vector superspace $H^0(\mathcal{G}_0, \mathcal{T})^{\mathcal{G}_0}$.
This action we can describe in terms of
the corresponding super Harish-Chandra pair $(\mathcal{G}_0,\mathfrak{g})$ in the following way:
\begin{equation}\label{eq action on left inv vectors}
g \mapsto (X \mapsto r^*_g \circ X \circ (r^{-1}_g)^*), \quad Y
\mapsto (X \mapsto [Y,X]),
\end{equation}
where $g\in \mathcal{G}_0$, $X\in H^0(\mathcal{G}_0, \mathcal{T})^{\mathcal{G}_0}$ and $Y\in \mathfrak{g}$.
Note that this action is well-defined because $\mathcal{G}$-left and right actions on $H^0(\mathcal{G}_0, \mathcal{T})$
commute. The Lie supergroup $\mathcal{G}$ acts also on
the vector superspace $H^0(\pt,\mathcal{O}_{\mathcal{G}_0\backslash
\mathcal{G}}) \otimes \mathfrak{g}$. This action is given by right translations $r^*_g$ on $H^0(\pt,\mathcal{O}_{\mathcal{G}_0\backslash
\mathcal{G}})$ and by the formulas (\ref{eq action on left inv vectors}) on $\mathfrak{g}$ if we assume that $X\in \mathfrak{g}$.
Clearly, the isomorphism $F$ from Lemma \ref{left invariant grading operators} is equivariant.
From now on we will identify  $H^0(\mathcal{G}_0, \mathcal{T})^{\mathcal{G}_0}$ and  $H^0(\pt, \mathcal{O}_{\mathcal{G}_0\backslash
\mathcal{G}}) \otimes \mathfrak{g}$ via isomorphism $F$ from Lemma \ref{left invariant grading operators}.

 If $\mathcal{H}$ is a Lie
subsupergroup of $\mathcal{G}$ and $\mathfrak{h} = \Lie \mathcal{H}$
then $\mathfrak{g}/\mathfrak{h}$ is an $\mathcal{H}$-module.

\medskip

\lem\label{lem w is H-invariant} {\it Let us take a $\mathcal{G}_0$-left invariant split grading operator $w$. The vector field $w$
satisfies (\ref{conditions v invariant}) if and only if
\begin{equation}\label{description of w in invariants}
\overline{w}\in \big(\!H^0(\pt, \mathcal{O}_{\mathcal{G}_0\backslash
\mathcal{G}})\otimes \mathfrak{g}/\mathfrak{h}\big)^{\mathcal{H}},
\end{equation}
 where
$\overline{w}$ is the image of $w$ by the natural mapping
$$
H^0\big(\!\pt, \mathcal{O}_{\mathcal{G}_0\backslash \mathcal{G}}\big)\otimes
\mathfrak{g} \to H^0\big(\!\pt, \mathcal{O}_{\mathcal{G}_0\backslash
\mathcal{G}}\big)\otimes \mathfrak{g}/\mathfrak{h}.
$$

 }

\medskip

\noindent{\it Proof.} Let $\overline{w}\in (H^0(\pt,
\mathcal{O}_{\mathcal{G}_0\backslash \mathcal{G}})\otimes
\mathfrak{g}/\mathfrak{h})^{\mathcal{H}}$. It follows that
$$
r_h^* \circ w \circ (r^{-1}_h)^* - w\in H^0\big(\!\pt,
\mathcal{O}_{\mathcal{G}_0\backslash \mathcal{G}}\big)\otimes
\mathfrak{h},\,\, h\in \mathcal{H}_0,
$$
and
$$
[Y,w] \in H^0\big(\!\pt, \mathcal{O}_{\mathcal{G}_0\backslash
\mathcal{G}}\big)\otimes \mathfrak{h},\,\, Y\in \mathfrak{h}.
$$
Hence, the conditions (\ref{conditions v invariant}) are satisfied.

On the other hand, if the conditions (\ref{conditions v invariant})
are satisfied, then the vector fields $r_h^* \circ w \circ
(r^{-1}_h)^* - w$ and $[Y,w]$ are vertical with respect to the
projection $\pi:\mathcal{G} \to \mathcal{G}/\mathcal{H}$. Therefore,
$r_h^* \circ w \circ (r^{-1}_h)^* - w$ and $[Y,w]$ belong to the superspace $H^0(\pt, \mathcal{O}_{\mathcal{G}_0\backslash \mathcal{G}})\otimes
\mathfrak{h}$. It is equivalent to conditions (\ref{description of w in invariants}).$\Box$

\medskip

Now our aim is to describe the space $\big(\!H^0(\pt,
\mathcal{O}_{\mathcal{G}_0\backslash \mathcal{G}})\otimes
\mathfrak{g}/\mathfrak{h}\big)^{\mathcal{H}_0}$. We have seen in
Proposition \ref{Ex_supergroups are nin-split} that the superspace $H^0(\pt,
\mathcal{O}_{\mathcal{G}_0\backslash \mathcal{G}})$ is isomorphic to
$\bigwedge \mathfrak{g}_{\bar 1}^*$. Actually this isomorphism can be chosen in
$\mathcal{G}_0$-equivariant way. More precisely, we need the following lemma.

\medskip

\l\label{prop left} {\it {\bf a.} We have
$$
\begin{array}{rll}
H^0\big(\!\pt, \mathcal{O}_{\mathcal{G}_0\backslash \mathcal{G}}\big)\otimes
\mathfrak{g} &\simeq\,\, \bigwedge (\mathfrak{g}_{\bar 1}^*)\otimes
\mathfrak{g} & \text{as
$\mathcal{G}_0$-modules},\\
\rule{0pt}{5mm} H^0\big(\!\pt, \mathcal{O}_{\mathcal{G}_0\backslash
\mathcal{G}}\big)\otimes \mathfrak{g}/\mathfrak{h}&\simeq\,\, \bigwedge
(\mathfrak{g}_{\bar 1}^*)\otimes \mathfrak{g}/\mathfrak{h} &
\text{as $\mathcal{H}_0$-modules},
\end{array}
$$
where the action of $\mathcal{G}_0$ on $\bigwedge (\mathfrak{g}_{\bar 1}^*)$ is standard.

{\bf b.} There exists a $\mathcal{G}_0$-left and right invariant split grading operator on $\mathcal{G}$.}

\medskip

\noindent{\it Proof.} {\bf a.} We have to show that there exists an $\mathcal{G}_0$-equivariant isomorphism
$$
H^0(\pt,\mathcal{O}_{\mathcal{G}_0\backslash \mathcal{G}}) \stackrel{\beta}{\longrightarrow}  \bigwedge \mathfrak{g}_{\bar 1}^*.
$$
Then the map $\beta\otimes \id$ will provide the required  isomorphism of $\mathcal{G}_0$-modules.
Consider the Lie supergroup
$$
\mathcal{G}^3= (\mathcal{G}_0, \Hom_{\mathbb C}(\bigwedge
\mathfrak{g}_{\bar 1}, \mathcal{F}_{\mathcal{G}_0}))
$$
from Section $3.4$. It
follows from (\ref{left right translations}) that
$$
H^0(\mathcal{G}_0,\mathcal{O}_{\mathcal{G}_0 \backslash \mathcal{G}^3}) =  \Hom_{\mathbb C}(\bigwedge
\mathfrak{g}_{\bar 1}, \mathbb C) = (\bigwedge
\mathfrak{g}_{\bar 1})^*.
$$
Note that the action of $\mathcal{G}_0$ on $(\bigwedge
\mathfrak{g}_{\bar 1})^*$ by right translations in $\mathcal{G}^3$,
denoted by $(r'_g)^*$,  coincides with the standard action of
$\mathcal{G}_0$ on $(\bigwedge \mathfrak{g}_{\bar 1})^*$. Indeed,
let us take
$$
f\in H^0(\mathcal{G}_0,\mathcal{O}_{\mathcal{G}^3})^{\mathcal{G}_0} = (\bigwedge
\mathfrak{g}_{\bar 1})^*.
$$
 By (\ref{left right translations split}), we have:
$$
(r'_g)^*(f)(X)(e) = (r'_g)^*(f)(X)(h) = f\big(\!\Ad(g^{-1}) X \big) (hg) = f\big(\!\Ad(g^{-1}) X \big) (e).
$$
Here $g,h\in \mathcal{G}_0$, $X\in \bigwedge
\mathfrak{g}_{\bar 1}$ and $e\in  \mathcal{G}_0$ is the identity.
It remains to note that by Lemma \ref{lem Phi is left right invariant}, the map
$\Phi_{\mathfrak{g}}$ induces the equivariant isomorphism  between the superspaces of left invariants $H^0(\pt,
\mathcal{O}_{\mathcal{G}_0\backslash \mathcal{G}})$ and $(\bigwedge
\mathfrak{g}_{\bar 1})^*$.

\smallskip

{\bf b.} We need to show that in the vector space $$
\big(\!\bigwedge (\mathfrak{g}_{\bar 1}^*)\otimes \mathfrak{g}\big)^{\mathcal{G}_0} =
\big(\!\bigwedge (\mathfrak{g}_{\bar 1}^*)\otimes \mathfrak{g}_{\bar 0}\big)^{\mathcal{G}_0}
\oplus \big(\!\bigwedge (\mathfrak{g}_{\bar 1}^*)\otimes \mathfrak{g}_{\bar 1}\big)^{\mathcal{G}_0}
$$
there exists points corresponding to split grading operators. This
space always possesses a $\mathcal{G}_0$-invariant, precisely, the
identity operator $\id \in \mathfrak{g}_{\bar 1}^*\otimes
\mathfrak{g}_{\bar 1}$. The pre-image of $\beta^{-1}(\id) \in
H^0(\pt, \mathcal{O}_{\mathcal{G}_0\backslash \mathcal{G}})\otimes
\mathfrak{g}$ has the form $\sum \varepsilon^i X_i$ for some choice
of local coordinates such that $X_i(\varepsilon^j)(e) =
\delta_{i}^j$, see Example \ref{Ex_split grading operator for split
Lie sgroup}. We have seen that such vector fields correspond to
$\mathcal{G}_0$-left invariant split grading operators on
$\mathcal{G}$.$\Box$

\medskip

Denote by $\mathcal{T}_{\mathcal{G}}$ the tangent sheaf of a Lie group $\mathcal{G}$ and by $\overline{v}$ is the image of $v$ by the natural mapping
$$
H^0(\pt, \mathcal{O}_{\mathcal{G}_0\backslash \mathcal{G}})\otimes
\mathfrak{g} \to H^0(\pt, \mathcal{O}_{\mathcal{G}_0\backslash
\mathcal{G}})\otimes \mathfrak{g}/\mathfrak{h}.
$$

The result of our study is:

\medskip

\t\label{theor grading operators on M} {\it The following conditions
are equivalent:

\smallskip

{\bf a.} A homogeneous supermanifold $\mathcal{M}=\mathcal{G}/\mathcal{H}$ admits a $\mathcal{G}_0$-left
invariant split grading that is induced by a grading of $\mathcal{O}_{\mathcal{G}}$
and the inclusion $\mathcal{O}_{\mathcal{M}}\subset (\pi_0)_*( \mathcal{O}_{\mathcal{G}})$.

\smallskip

{\bf b.} There exists a $\mathcal{G}_0$-left invariant vector field  $\chi \in H^0(
\mathcal{G}_0,(\mathcal{T}_{\mathcal{G}})_{(2)\bar 0})$ such that
\begin{equation}\label{w chi 1}
\overline{\chi} \in  \left(\bigwedge
(\mathfrak{g}_{\bar 1})^* \otimes \mathfrak{g}/\mathfrak{h}
\right)^{\mathcal{H}_0},
\end{equation}
and such that for $w = \beta^{-1}(\id) + \chi$, where $\beta^{-1}(\id) =\sum \varepsilon^i
X_i$ is from the proof of Proposition \ref{prop left}.b, we have
\begin{equation}\label{w 2}
[Y,w] \in H^0(\pt, \mathcal{O}_{\mathcal{G}_0\backslash
\mathcal{G}}) \otimes \mathfrak{h}, \,\,\,Y\in \mathfrak{h}_{\bar
1}.
\end{equation}
 }

\section{\bf An application}

As above let $\mathcal{G}$ be a Lie supergroup and $\mathcal{H}$ be a Lie subsupergroup, $\mathfrak{g}$ and
$\mathfrak{h}$ be the Lie superalgebras of $\mathcal{G}$ and
$\mathcal{H}$, respectively, and $\mathcal{M}:=
\mathcal{G}/\mathcal{H}$. Consider the map
$$
\rho: \mathfrak{g}_{\bar
0} \to H^0(\pt, \mathcal{T}_{\mathcal{G}_0\backslash \mathcal{G}})
$$
induced by the action of $\mathcal{G}_0$ on $\mathcal{M}$. (Here
$\mathcal{T}_{\mathcal{G}_0\backslash \mathcal{G}}$ is the sheaf of
vector fields on $\mathcal{G}_0\backslash \mathcal{G}$.) Let us
describe its kernel. For $X\in \mathfrak{g}_{\bar 0}$ and
$f\in H^0(\pt, \mathcal{O}_{ \mathcal{G}_0\backslash \mathcal{G}})$, we have:
$$
\begin{array}{rl}
X(f)(Y)(e) = &\frac{\d}{\d t}|_{t=0}
f\big(\!\Ad(\exp(-tX)) Y\big)\big(\!\exp(tX)\big) =\\
\rule{0pt}{6mm}&\frac{\d}{\d t}|_{t=0} f\big(\!\Ad(\exp(-tX)) Y\big)(e),
\end{array}
$$
where $Y= Y_1\cdots Y_r$, $Y_i\in \mathfrak{g}_{\bar 1}$ and $t$ is
an even parameter. A vector field $X$ is in $\Ker \rho$ if and only if
$X(f)(Y)(e)= 0$ for all  $f$ and $Y$. Hence,
$$
\Ker \rho = \Ker (\ad|_{\mathfrak{g}_{\bar 1}}),
$$
where $\ad$ is the adjoint representation of $\mathfrak{g}_{\bar 0}$
in $\mathfrak{g}$.

Furthermore, denote
$$
\begin{array}{rl}
A:= &\Ker \big(\mathcal{G}_0\ni g \mapsto \overline{l}_g :
\mathcal{G}/\mathcal{H} \to \mathcal{G}/\mathcal{H}\big);\\
\rule{0pt}{5mm}\mathfrak{a}: = &\Ker \big(\mathfrak{g}\ni X \mapsto
H^0(\mathcal{G}_0/\mathcal{H}_0,
\mathcal{T}_{\mathcal{G}/\mathcal{H}})\big).
\end{array}
$$
Here $\overline{l}_g$ is the automorphism of
$\mathcal{G}/\mathcal{H}$ induced by the left translation $l_g$. The
pair $(A,\mathfrak{a})$ is a super Harish-Chandra pair. An action of
$\mathcal{G}$ on $\mathcal{M}$ is called {\it effective} if the
corresponding to $(A,\mathfrak{a})$ Lie supergroup is trivial. As in
the case of Lie groups any action of a Lie supergroup can be factored to be effective.

\medskip

\t\label{main split} {\it Assume
that the action of $\mathcal{G}$ on $\mathcal{M}$ is effective. If
$$
[\mathfrak{g}_{\bar 1}, \mathfrak{h}_{\bar 1}]\subset
\mathfrak{h}_{\bar 0} \cap \Ker (\ad|_{\mathfrak{g}_{\bar 1}}),
$$
then $\mathcal{M}$ is split.}

\medskip

\noindent{\it Proof.} Let us show that in this case the vector field
$w =\sum \varepsilon^iX_i +0 = \sum \varepsilon^iX_i$ from Proposition \ref{prop left}.b is a (left invariant) split grading operator on
$\mathcal{M}$ using Theorem \ref{theor grading operators on M}.

The
condition (\ref{w chi 1}) is satisfied trivially, because $\chi =0$.
Let us check the condition (\ref{w 2}). We have:
$$
[Y,v] = \sum Y(\varepsilon^i) X_i - \sum \varepsilon^i [Y,X_i].
$$
Since $[\mathfrak{g}_{\bar 1}, \mathfrak{h}_{\bar 1}]\subset
\mathfrak{h}_{\bar 0}$, we get
$$
\sum \varepsilon^i [Y,X_i]\in
H^0(\pt, \mathcal{O}_{\mathcal{G}_0\backslash \mathcal{G}}) \otimes
\mathfrak{h}.
$$
 Hence, we have to show that
 $$
 \sum Y(\varepsilon^i)
X_i\in H^0(\pt, \mathcal{O}_{\mathcal{G}_0\backslash \mathcal{G}})
\otimes \mathfrak{h}.
$$

Assume that $X_1,\dots X_k$ is a basis of $ \mathfrak{h}_{\bar 1}$, $X_1,\dots X_k, X_{k+1},\dots X_m$
is a basis of $\mathfrak{g}_{\bar 1}$ and $(\varepsilon^i)$ is the system of global odd $\mathcal{G}_0$-left
invariant coordinates corresponding to this basis such that $\sum \varepsilon^iX_i$ is as in Proposition
\ref{prop left}.b. In particular, $\varepsilon^i (\gamma_{\mathfrak{g}}( X_j))= \delta^i_j$, because
$\sum (\varepsilon^i\circ \gamma_{\mathfrak{g}})\otimes X_i$ is the identity operator in
$\mathfrak{g}_{\bar 1}^* \otimes \mathfrak{g}_{\bar 1}$.

Let us take $Z\in \Ker \rho$. Clearly, $Z(\varepsilon^i) =0$ and $X_j (\varepsilon^i)$ is
again a $\mathcal{G}_0$-left invariant function on $\mathcal{G}$. By (\ref{left inv vector field}), we also have:
$$
\varepsilon^i(X_{i_1} \cdots Z \cdots  X_{i_k}) =0.
$$
Furthermore, by definition of $\varepsilon^i$, we get that
$\varepsilon^i\circ \gamma_{\mathfrak{g}}\in \mathfrak{g}_{\bar
1}^*$. Hence,
$$
\varepsilon^i(\gamma_{\mathfrak{g}}( X_{i_1}\wedge  \cdots \wedge X_{i_k})) = 0,
$$
if $k>1$. Summing up all these observations we see that
$$
\varepsilon^i(\gamma_{\mathfrak{g}}( X) \cdot Y) = \varepsilon^i(\gamma_{\mathfrak{g}}( X \wedge Y)) +0,
$$
where $Y\in \mathfrak{h}$ and $X\in \bigwedge \mathfrak{g}_{\bar 1}$. Now we can conclude that

$$
\sum Y(\varepsilon^i) X_i =  - Y \in \mathfrak{h}\subset H^0(\pt,
\mathcal{O}_{\mathcal{G}_0\backslash \mathcal{G}}) \otimes
\mathfrak{h}.
$$
The proof is complete.$\Box$

\medskip
\smallskip

\ex \label{ex split super-grassmanians} Consider the
super-grassmannian $\mathbf{Gr}_{m|n,k|l}$. It is a
$\mathrm{GL}_{m|n}$-homogeneous space, see \cite{onigl} for more
details. Hence, $\mathbf{Gr}_{m|n,k|l} \simeq
\mathrm{GL}_{m|n}/\mathcal{H}$ for a certain $\mathcal{H}$. (See,
for example, \cite{V_Func}.) It the case $k=0$ or $k=m$, the
following holds $[(\mathfrak{gl}_{m|n})_{\bar 1}, \mathfrak{h}_{\bar
1}]=0$. Therefore, by Theorem \ref{main split}, the
super-grassmannian is split.

In \cite{onigl} it was shown that the super-grassmannian
$\mathrm{GL}_{m|n,k|l}$ is not split if and only if $0<k<m$ and
$0<l<n$. (This fact also follows from results in \cite{LPW} and
\cite{PS} about non-projectivity of super-grassmannian.)

\medskip

Finally, let us recall a result proved in \cite{V_TG}:

\medskip

\t\label{TG} {\it If a complex homogeneous supermanifold
$\mathcal{M}$ is split, then there is a Lie supergroup $\mathcal{G}$
with $[\mathfrak{g}_{\bar 1}, \mathfrak{g}_{\bar 1}] = 0$, where
$\mathfrak{g}=\mathfrak{g}_{\bar 0}\oplus \mathfrak{g}_{\bar 1} =
\Lie \mathcal{G}$, such that $\mathcal{G}$ acts on $\mathcal{M}$
transitively. }

\noindent{\it Elizaveta Vishnyakova}

\noindent {Max Planck Institute for Mathematics Bonn and}

\noindent University of Luxembourg

 \noindent {\emph{E-mail address:}
\verb"VishnyakovaE@googlemail.com"}

\end{document}